\documentclass[11pt]{article}
\usepackage{amssymb}
\usepackage{amsbsy}
\begin{document}
\centerline{\Large\bf The John Theorem for Simplex
\footnote{Supported by
National Natural Sciences Foundation of China(~No.10271071) and
the Youth Science Foundation of Shanghai(~No.214511)}} \vskip 20pt
\begin{center}

{ Si Lin$^{1,2}$, ~~Xiong Ge$^2$ and ~Leng Gangsong$^2$}

\vskip 16pt

{\it 

1. Department of Mathematics, Beijing Forestry University, Beijing,100083, China

2. Department of Mathematics, Shanghai University, Shanghai,200444, China}

\vskip 16pt
\begin{minipage}{11cm}
{\bf Abstract}.~ { In this paper, we give a description of the
John contact points of a regular simplex. We prove that the John
ellipsoid of any simplex is ball if and only if this simplex is
regular and that the John ellipsoid of a regular simplex is its
inscribed ball.}



\end{minipage}
\end{center}

\vskip 20pt {\center\large \bf 1.~~Introduction }

\vskip 10pt

\vskip 5pt

In 1948, F. John proved that every convex body (i.e., a compact,
convex subset with nonempty interior) in $R^n$ contains only one
maximal(in volume) ellipsoid, which is known as the {\it John
ellipsoid}. When the John ellipsoid is the unit ball $B_2^n$,
F.John has proved the following theorem.

\noindent{\bf Theorem 1.[John]}~~{\it Let $C$ be a convex body in
$R^n$. The ellipsoid of maximal volume in $C$ is $B^n_2$, if and
only if $C$ contains $B_2^n$ and there are some points $(u_i)_1^m$
on the boundary of $C$ and positive numbers $(c_i)_1^m$ so that}
$$\begin{array}{rl}
&a)~~\displaystyle \sum_{i=1}^{m}c_iu_i\otimes u_i=I_n,~~~ and \\
&b)~~\displaystyle \sum_{i=1}^{m}c_iu_i=0.
\end{array}$$

\noindent Here, $I_n$ is the identity map on $R^n$ and, for any
unit vector $u$, $u\otimes u$ is the rank-one orthogonal
projection onto the span of $u$, i.e., the map $x\longrightarrow
\langle x, u\rangle u.$ The $u_i^,$s of the theorem are the
intersection points of the unit sphere $S^{n-1}$ with the boundary
of $C$.

Condition a) shows that the $(u_i)_1^m$ behave rather like an
orthonormal basis in that we can resolve the Euclidean norm as a
(weighted) sum of square of inner products. This condition is
equivalent to the statement that, for all $x$ in $R^n$,

$$|x|^2=\sum_{i=1}^m c_i\langle x,u_i\rangle^2 ,\eqno(1)$$
where $\langle,\rangle$ is the usual Euclidean inner product and
$|\cdot|$ is the induced norm by this inner product.

By a simple computation, we know that equality (1) is equivalent
to
$$x=\sum_{i=1}^m c_i\langle x,u_i\rangle u_i. \eqno(2)$$

For detail, see [Ba1] and in [Ba3], one can find a modern proof of
Theorem 1.

\noindent{\bf Definition 1}~~{\it Suppose that
$A_1,...,A_{n+1}~\in R^n$ be affinely independent then the convex
hull, denoted by $A$, of these points is called a simplex, i.e.,
$${A}=\{x|~x=\sum_{i=1}^{n+1}\lambda_i A_i, ~~\sum_{i=1}^{n+1}\lambda_i=1, ~\lambda_i\geq 0\},$$
and if all $|A_iA_j|, ~i\neq j$ are equal, then we call $A$ a
regular simplex. }

In Theorem 1., $(u_i)_1^m$ is usually called the {\it contact
points}. For the unit cube $[-1,1]^n$ in $R^n$, the maximal
ellipsoid is $B_2^n$ as one would expect, so the contact points
are the standard basis vectors $(e_1,...,e_n)$ of $R^n$ and their
negatives. However, even for the simplest nonsymmetric convex
body---simplex, there is no nature description of the contact
points.

In this paper, we give a description of the contact points for a
regular simplex and the main results are the following theorems.

\noindent{\bf Theorem 2.}~~{\it The John ellipsoid of a regular
simplex is its inscribed ball.}

\noindent{\bf Theorem 3.}~~{\it For any simplex in $R^n$, the John
ellipsoid of this simplex is ball if and only if the simplex is
regular.}

\vskip 30pt

{\center\large \bf 2.~~The Proof of Main Results }

\vskip 5pt

First we introduce the following definition[C].

\noindent{\bf Definition 2[C]}~~{\it Suppose that $A$ is an
n-dimensional simplex with vertexes $\{A_1,...A_{n+1}\}$, $M$ is a
point in $R^n$. Denote by $V_i, i=1,...,n+1$, the volume of the
simplex with vertexes $\{A_1,...A_{i-1},M$ $,
A_{i+1},...,A_{n+1}\}$, and if the dimensions of
$$con{\{A_1,...A_{i-1},M, A_{i+1},...,A_{n+1}\}} ~~{ and}~~
con{\{A_1,...A_{i-1},M, A_{i+1},...,A_{n+1}\}}\cap A$$ are both n,
then the following ratio is called the {\it barycentric
coordinates} of $M$,}
$$V_1:V_2:\ldots:V_{n+1}.$$

Suppose that $\{A_1,...A_{n+1}\}$ are the vertexes of a regular
simplex $A$ and $B_2^n$ is its inscribed ball. Denote by
$\{B_i,i=1,...,n+1\}$, the tangent points of $B_2^n$ with the face
generated by the convex hull of $
\{A_1,...,A_{i-1},A_{i+1},...,A_{n+1}\}$ and denote by
$\{u_i,i=1,...,n+1\}$ the outer normal unit vectors of these facet
respectively. According to the Definition 2, we have the
barycentric coordinates of $B_i$ as follows
$$(\underbrace{\frac{1}{n},...,\frac{1}{n},0,\frac{1}{n},...,\frac{1}{n}}_{n+1}),$$
where $0$ is in the $i$-th$(i=1,...,n+1)$ position.

\vskip 5pt

\vskip 12pt

\vskip 10pt \noindent {\bf Proof of Theorem 2.}

\vskip 5pt

According to the Theorem 1., it suffices to prove that the tangent
points of a regular simplex with its inscribed ball satisfied the
condition a) and b).

Now suppose that $A$ is a regular simplex with vertexes
$\{A_1,A_2,...,A_{n+1}\}$ and $B_2^n$ is its inscribed ball.
Denote by $\{B_i,i=1,...,n+1\}$ the tangent points which is
opposite to $\{A_i,i=1,...,n+1\}$ respectively. From the above
discussion, the barycentric coordinates of $B_i$ is
$$(\underbrace{\frac{1}{n},...,\frac{1}{n},0,\frac{1}{n},...,\frac{1}{n}}_{n+1}),$$
where $0$ is in the $i$-th$(i=1,...,n+1)$ position.

Obviously, the barycentric coordinates of the origin is
$(\underbrace{1,...,1}_{n+1}).$

Let $c_i=\frac{n}{n+1},i=1,...,n+1$, then
$$\sum_{i=1}^{n+1}c_iB_i=
\underbrace{(1,...,1)}_{n+1}.$$ Thus the condition b) is
satisfied.

Next, we will prove that $u_i$ satisfied the condition a) of
Theorem 1., that is, for any $x\in R^n$, the following equality
holds,

$$x=\sum_{i=1}^{n+1}c_i\langle x,
u_i\rangle u_i,$$
where $c_i=\frac{n}{n+1},i=1,...,n+1.$

Because $A$ is a $n$-dimensional simplex, the space spaned by the
$n+1$ vectors $\{u_i,i=1,...,n+1\}$ must be $R^n$, i.e.,
$$Span\{u_1,u_2,...,u_{n+1}\}=R^n.$$
So for any $x\in R^n$, there must exist $n+1$ real numbers
$\alpha_1,....,\alpha_{n+1}$ such that
$$x=\alpha_1u_1+....+\alpha_{n+1}u_{n+1}.\eqno(3)$$

Thus, we can get
\begin{eqnarray}
\left\{
\begin{array}{l}
\langle u_1,x\rangle=\alpha_1\langle u_1,u_1\rangle+\alpha_2\langle u_1,u_2\rangle+....+\alpha_{n+1}\langle u_1,u_{n+1}\rangle, \\
\langle u_2,x\rangle=\alpha_1\langle u_2,u_1\rangle+\alpha_2\langle u_2,u_2\rangle+....+\alpha_{n+1}\langle u_2,u_{n+1}\rangle, \\
~~~~~\vdots\\

\langle u_{n+1},x\rangle=\alpha_1\langle u_{n+1},u_1\rangle+\alpha_2\langle u_{n+1},u_2\rangle+....+\alpha_{n+1}\langle u_{n+1},u_{n+1}\rangle. \\
\end{array}
\right.  \label{3}\nonumber
\end{eqnarray}

\noindent Denote $\alpha=(\alpha_1,....,\alpha_{n+1}),
\beta=(\langle u_1,x\rangle,....,\langle u_{n+1},x\rangle),$ and

$$D=\left(\begin{array}{cccccccc}
{\langle u_1,u_1\rangle} & {\displaystyle \langle u_1,u_2\rangle}&
\cdots &
{\displaystyle \langle u_1,u_{n+1}\rangle}\\

{\langle u_2,u_1\rangle} & {\displaystyle \langle u_2,u_2\rangle}&
\cdots &
{\displaystyle \langle u_2,u_{n+1}\rangle}\\

\cdots&\cdots&\cdots&\cdots\\

{\langle u_{n+1},u_1\rangle} & {\displaystyle \langle
u_{n+1},u_2\rangle}& \cdots & {\displaystyle \langle
u_{n+1},u_{n+1}\rangle}

\end{array}\right)_{(n+1)\times(n+1)},$$

\noindent then the above equation system can be written as
$$D\alpha^T=\beta^T,\eqno(4)$$
\noindent where $\alpha^T,\beta^T$ represent respectively the
transform of $\alpha$ and $\beta$.

Observe that every element of $D$, $\langle u_{i},u_j\rangle$, is
the cosine of angle of two outer normal unit vectors. Denote by
$F_i,F_j$ the faces whose outer normal unit vectors are
$u_{i},u_j$ respectively. Obviously, the angle of $u_{i},u_j$ is
mutually complementary with the angle of $F_i,F_j$, i.e.,
$$\langle u_{i},u_j\rangle=-\cos\angle(F_i,F_j),$$
where $\angle(F_i,F_j)$ represents the dihedral angle of
$F_i,F_j$.

For $\cos\angle(F_i,F_j)$, we have the following equality,
$$\cos\angle(F_i,F_j)=\frac{S_{ji}}{S_j},$$
where $S_j$ is the $(n-1)$-dimensional volume of face $F_j$, and
$S_{ji}$ is the volume of the projection $F_j$ to $F_i$ along
$u_i$ .

For $A$ is regular simplex, the $n ~(n-1)$-dimensional volumes of
all the projections $F_j, j\neq i$ to $F_i$ along $u_i$ are equal.
So $\frac{S_{ji}}{S_j}=\frac{1}{n}.$ Thus we get the $D$, i.e.

$$D=\left(\begin{array}{cccccccc}
{1} & {\displaystyle -\frac{1}{n}}& \cdots &
{\displaystyle -\frac{1}{n}}\\

{\displaystyle -\frac{1}{n}}&1&\cdots&
{\displaystyle -\frac{1}{n}}\\

\cdots&\cdots&\cdots&\cdots\\

{\displaystyle -\frac{1}{n}}& {\displaystyle -\frac{1}{n}}&
\cdots&1

\end{array}\right)_{(n+1)\times(n+1)}.$$

It follows from condition b) and (4) that
\begin{eqnarray}
\left\{
\begin{array}{l}
D\alpha^T=\beta^T\\
\sum_{i=1}^{n+1}\langle u_i,x\rangle=0.\\
\end{array}
\right.  \label{21}\nonumber
\end{eqnarray}

Let $\alpha=(\frac{n}{n+1}\langle
u_1,x\rangle,...,\frac{n}{n+1}\langle u_{n+1},x\rangle)$ in the
above equation system, we know that $\alpha$ is a solution of this
equation system. So every point $x \in R^n$ can be represented as
the form of (2).

\noindent The proof of Theorem 2. is completed.

\vskip 10pt

To prove Theorem 3., we need the following Brascamp-Lieb
inequality, which is the generalization of covolution inequality.

\noindent{\bf Theorem 4.[BL]}~~Suppose that $(u_i)^m_1$ is a
sequence of unit vector in $R^n$, $(c_i)^m_1$ is a sequence
positive real numbers and they satisfied the following equality
$$\sum_{i=1}^{m}c_iu_i\otimes u_i=I_n.$$
If $f_i: R\longrightarrow [0,\infty), i=1,...,m$ is a sequence of
integrable functions, then
$$\int_{R^n}\prod_{i=1}^mf_i(\langle u_i, x\rangle)^{c_i}dx\leq \prod_{i=1}^m(\int_Rf_i)^{c_i}. \eqno(5)$$

F.Barthe get a necessary condition for the equality holds in
Theorem 4.

\noindent{\bf Theorem 5.[Bar]}~~Suppose that $(u_i)^m_1$ is a
sequence of unit vector in $R^n$, $(c_i)^m_1$ a sequence positive
real numbers, and they satisfied the following equality
$$\sum_{i=1}^{m}c_iu_i\otimes u_i=I_n.$$
If $(f_i)_1^m$ is a sequence functions, not all zero in $L_1(R)$,
and all $(f_i)_1^m$ are not the density function of Gauss
distribution, then the necessary condition for the equality hold
in (5) is
$$m=n,$$
and $$(u_i)_1^m$$ is a orthonormal basis of $R^n$.

\vskip 30pt

\noindent {\bf Proof of Theorem 3.}

\vskip 5pt

The "if" part of Theorem 3. can be obtained from Theorem 2.
directly. So it is sufficient to prove that if the John ball of
the simplex $C$ is $B_2^n$ then $C$ is regular.

Firstly, we observe that if the John ball of the simplex $C$ is
$B_2^n$, then $B_2^n$ is the inscribed ball of $C$. If not,
without lost of generalization, suppose that $B_2^n$ is not
tangent with face $F_i$. Let $u_i$ be the outer normal unit vector
of $F_i$, then there must exist a positive number $\varepsilon$,
such that $B_2^n$ is not tangent with any faces of $C$ when
$B_2^n$ move $\varepsilon$ along $u_i$. At this time, there must
exist another positive number $r>1$ such that $rB_2^n$ be the John
ball of $C$. This contradicts with the fact that $B_2^n$ is the
John ball of the simplex $C$.

Because the inscribed ball of $C$ is it's John ball, by Theorem
1., there exist a sequence positive real numbers $(c_i)_1^{n+1}$
and a sequence of unite vectors $(u_i)_1^{n+1}$ on the boundary of
$C$ such that
$$\sum_{i=1}^{n+1}c_iu_i\otimes u_i=I_n, \eqno(6)$$
and
$$\sum_{i=1}^{n+1}c_iu_i=0.\eqno(7)$$
Denote by $K=\{x\in R^n: \langle x, u_i\rangle\leq1, 1\leq i\leq
n+1\}$, then $K$ is also the simplex in $R^n$. Because
$(u_i)_1^{n+1}$ are the contact points of $C$ and $B_2^n$,
$$C\subset \{x\in R^n: \langle x, u_i\rangle\leq1, 1\leq i\leq m\}=K.$$
Observe that $B^n_2$ is also the inscribed ball of $C$ and that
$K$,$C$ have the same tangent points $(u_i)_1^{n+1}$ with $B_2^n$,
so
$$C=K.$$

Next, we will show that $K$ is regular simplex.

In the following discussion, $R^{n+1}$ will be regarded as
$R^n\times R.$ For each $i$ let
$$v_i=\sqrt{\frac{n}{n+1}}(-u_i, \frac{1}{\sqrt{n}})\in R^{n+1}, ~~~i=1,...,n+1,$$
$$d_i=\frac{n+1}{n}c_i, ~~~i=1,...,n+1,$$
then $v_i$ is a unit vector and the identities (6) and (7),
together yield that
$$\sum_{i=1}^{n+1}d_iv_i\otimes v_i=I_{n+1}.$$
Define a sequence functions $(f_i)_1^{n+1}$ as follows,
\begin{eqnarray}
f_i(t)=\left\{
\begin{array}{l}
e^{-t}, ~~~if ~~t\geq 0, \\
0, ~~~if ~~t< 0. \\
\end{array}
\right.  \nonumber
\end{eqnarray}
For any $x\in R^{n+1}$, let
$$F(x)=\prod_{i=1}^{n+1}f_i(\langle v_i, x\rangle)^{d_i},$$
by Theorem 4., we have
$$\int_{R^n}F(x)dx\leq \prod_{i=1}^{n+1}(\int_Rf_i)^{d_i}=1.\eqno(8)$$
Some of the above technique are from Ball. Using the similar
discussion in [Ba2], we get the integration of $F$ in the
hyperplane $\{x: x_{n+1}=r\geq0\}$
$$e^{-\sqrt{n+1}r}Vol(\frac{r}{\sqrt{n}}K)=e^{-\sqrt{n+1}r}(\frac{r}{\sqrt{n}})^nVol(K).$$
So by (8)
$$1\geq Vol(K)\int_0^{\infty}e^{-\sqrt{n+1}r}(\frac{r}{\sqrt{n}})^ndr=\frac{Vol(K)n!}{\sqrt{n^n(n+1)^{n+1}}},$$
i.e.,
$$Vol(K)\leq \frac{\sqrt{n^n(n+1)^{n+1}}}{n!}.\eqno(9)$$
Observe that the right hand of (9) is exactly the volume of the
regular simplex whose inscribed ball is $B_2^n$.

Observe the construction of $(f_i)_1^{n+1}$, and the Theorem 5.
for (8),thus we can get the condition for the equality holds in
(9) and that is $(v_i)_1^{n+1}$ is a sequence of orthonormal basis
of $R^{n+1}$. For any two vectors of this basis
$$v_i=\sqrt{\frac{n}{n+1}}(-u_i, \frac{1}{\sqrt{n}}),$$
and
$$v_j=\sqrt{\frac{n}{n+1}}(-u_j, \frac{1}{\sqrt{n}}),$$
we have
$$0=\langle v_i, v_j\rangle=\frac{n}{n+1}\langle u_i, u_j\rangle+\frac{1}{n}.$$
So
$$\langle u_i, u_j\rangle=-\frac{n+1}{n^2}, i\neq j,$$ is a
constant. Because that
 $(u_i)_1^{n+1}$ are the normal vectors of the $n+1$ faces of the simplex $K$, $K$ is a regular simplex.

\noindent This completes the proof of the theorems.

{\bf Acknowledgment}{~~The authors thank Dr. He Binwu for his
valuable advice for this paper.}

\vskip 10pt

\end{document}